\documentclass[10pt,a4paper]{amsart}
\usepackage[utf8]{inputenc}
\usepackage[T1]{fontenc}
\usepackage{amsmath}
\usepackage{amsfonts}
\usepackage{amssymb}
\usepackage{graphicx}
\usepackage{mathrsfs}
\usepackage{color}
\usepackage{tikz}

\usepackage{accents}
\newlength{\dhatheight}
\newcommand{\doublehat}[1]{%
	\settoheight{\dhatheight}{\ensuremath{\hat{#1}}}%
	\addtolength{\dhatheight}{-0.35ex}%
	\hat{\vphantom{\rule{1pt}{\dhatheight}}%
		\smash{\hat{#1}}}}

\newtheorem{theorem}{Theorem}[section]
\newtheorem{lemma}[theorem]{Lemma}
\newtheorem{corollary}[theorem]{Corollary}
\newtheorem{proposition}[theorem]{Proposition}

\theoremstyle{definition}

\newtheorem{definition}[theorem]{Definition}
\newtheorem{example}[theorem]{Example}

\theoremstyle{remark}

\newtheorem{Claim}{Claim}
\newtheorem{remark}[theorem]{Remark}

\DeclareMathOperator{\Aut}{{\rm Aut}}

\newcommand{\con}[1]{\mathop{\sf con}({#1})}
\newcommand{\prb}[1]{\mathop{\sf par}({#1})}

\newcommand{\lev}[1]{\mathop{\sf lev}({#1})}

\newcommand{\nub}[1]{\mathop{\sf nub}({#1})}
\newcommand{\lnub}[1]{\mathop{\sf lnub}({#1})}

\newcommand{\tdlc}{t.d.l.c.}

\newcommand{\frk}[1]{\mbox{\rm f-rk}(#1)}

\def\ident{1}
\def\triv{\{\ident\}}
\def\tree{T}

\author{George A. Willis}
\title{Groups with flat-rank greater than 1}

\thanks{This research was supported by the Australian Research Council grant FL170100032.}

\begin{document}
\begin{abstract}
A general method for finding subgroups of a totally disconnected, locally compact groups having flat-rank greater than 1 is described. This method uses the internal structure of the group, notably the Levi subgroup of a given flat group, in order to produce a new flat group that, under an additional hypothesis, has flat-rank one greater than that of the given group. 
\end{abstract}

	\maketitle

\section{Introduction}

A subgroup, $H$, of a totally disconnected, locally compact (\tdlc) group $G$ is \emph{flat} if there is a compact open subgroup, $U\leq G$ that is minimising for every $x\in H$, that is,  $U$ realises the minimum value of the set of positive integers 
$$
\left\{[xVx^{-1} : xVx^{-1}\cap V] \mid V\leq G\text{ compact and open}\right\}.
$$ 
This minimum value is called the \emph{scale of $x$} and denoted $s_G(x)$. The \emph{uniscalar subgroup} of $H$ is $H_u := \left\{h\in H\mid hUh^{-1} = U\right\}$ and $H/H_u$ is a free abelian group, \cite[Corollary 6.15]{SimulTriang}. The \emph{flat-rank} of $H$ is the rank of this free abelian group. Every singly generated group, $\langle x\rangle$, is flat with flat-rank equal to $0$ if there is $U\leq G$ normalised by $x$ and equal to $1$ otherwise. This paper is concerned with finding subgroups of $G$ having flat-rank greater than $1$. 

Known examples of flat groups having flat-rank greater than $1$ occur in cases where the ambient group $G$ has additional structure with which the flat group is associated. For example, if $G$ is $\mathrm{PSL}_n(\mathbf Q_p)$, the group of diagonal matrices is flat and has flat-rank $n-1$, equal to the Lie rank of $\mathrm{PSL}_n(\mathbf Q_p)$, see \cite[Example 6.11]{SimulTriang}. That is a special case of a general result about semisimple algebraic groups over local fields that is itself a consequence of the fact that, if $G$ is the automorphism group of a building and the building has an apartment with geometric rank $r$, then $G$ has abelian subgroup with flat-rank $r$, see \cite[Theorems A \& B]{BaumRemyWill} and \cite[Theorem~1]{Cap_Hag}. These results, showing that the flat-rank coincides with familiar notions of rank in important classes of groups, motivate the study of flat subgroups and their flat-rank in general \tdlc~groups.

There is no method for finding flat groups of higher rank in general \tdlc~groups however. In the case where $G$ is a Lie group over a local field and $x\in G$, $\mathrm{Ad}(x)$ is a linear operator on the Lie algebra of $G$ and the algebra of operators generated by $\mathrm{Ad}(x)$ is commutative. This commutative algebra is the Lie algebra of an abelian subgroup of $G$ whose flat-rank may be greater than $1$ and, moreover, is equal to the number of absolute values of eigenvalues of $\mathrm{Ad}(x)$, which may be found with the aid of the Cayley-Hamilton Theorem. While this method does not extend to general \tdlc~groups, it is the case that if the centraliser of $x$, $C(x)$, contains $y$ with $s_{C(x)}(y)>1$, then the subgroup $\langle x,y\rangle$ is flat and has flat-rank equal to that of $\langle x\rangle$ augmented by $1$. The main result shown below, Theorem~\ref{thm:extend_flat}, improves on this idea by replacing the centraliser of $x$ by the larger \emph{Levi subgroup} for $x$, that is, the set of all $y$ such that $\left\{x^nyx^{-n}\right\}_{n\in\mathbf{Z}}$ is precompact. It thus relates flat groups to Levi and parabolic subgroups in $G$, which were defined in \cite{ContractionB} in analogy with the similarly named subgroups of algebraic groups. The flat group found contains $y$ and a modification of $x$. 

The article has the following structure. Basic results about the scale and flat groups are recalled in \S\ref{sec:scale} and in a couple of cases sharpened. Preliminary results linking Levi subgroups and minimising subgroups are shown in \S\ref{sec:Levi} and \S\ref{sec:Main} gives the proof of the main theorem. A symmetrisation of the relation `$y$ is in the Levi subgroup of $x$' is suggested by the proof of the main theorem, and this symmetrised relation is explored in \S\ref{sec:symmetric_relation}. 
\section{Scale techniques}
\label{sec:scale}

This section collects results about the scale of an automorphism and related ideas that are used in the proof of the main theorem. The terms used here have in some cases only been introduced since the papers cited were written. Lemma~\ref{lem:tidy_criterion} strengthens known results and Theorem~\ref{thm:is_flat} is new (albeit fairly obvious).

\subsection{Tidy subgroups}
\label{sec:tidy}
A structural characterisation of subgroups minimising for an automorphism is given in \cite{FurtherTidy} in terms of tidiness of the subgroup, and tidiness of $U$ for $\alpha$ is defined in \cite{Structure94} in terms of the subgroups $U_+ = \bigcap_{k\in\mathbf{N}} \alpha^k(U)$ and $U_- := \bigcap_{k\in\mathbf{N}} \alpha^{-k}(U)$ as follows. (Zero is a natural number in this definition.)
\begin{definition}
	Let $\alpha\in\Aut(G)$ and $U$ be a compact open subgroup of $G$. Then $U$ is \emph{tidy for $\alpha$} if and only if
	\begin{description}
		\item[TA] $U = U_+U_-$ and
		\label{defn:TA}
		\item[TB] $U_{++} := \bigcup_{k\in\mathbf{N}} \alpha^k(U_+)$ is closed. 
			\label{defn:TB}
	\end{description}
\end{definition}

Then \cite[Theorem~3.1]{FurtherTidy} says that $U\leq G$ is minimising for $\alpha$ if and only if it is tidy for $\alpha$. Moreover, $\alpha(U_+)\geq U_+$ and $s(\alpha) = [\alpha(U_+) : U_+]$ if $U$ is tidy for $\alpha$. 

The compact open subgroup $U$ is said to be \emph{tidy above} for $\alpha$ if it satisfies {\bf TA} and \emph{tidy below} if it satisfies {\bf TB}.
The following is \cite[Lemma 1]{Structure94}. 
\begin{lemma}
	\label{lem:tidy_above}
Let $U$ be a compact open subgroup of $G$ and $\alpha\in \Aut(G)$. Then there is $n\geq0$ such that $\bigcap_{k=0}^n \alpha^k(U)$ is tidy above for $\alpha$.
\end{lemma}

It follows immediately from the definition that, if $U$ is minimising for $\alpha$, then so is $\alpha(U)$, and the same thus holds for subgroups tidy for $\alpha$. On the other hand, it may be shown, see \cite[Lemma 10]{Structure94}, that the intersection of two tidy subgroups for $\alpha$ is tidy for $\alpha$. The same thus holds for minimising subgroups and this fact will be used in what follows. Given $U$ be tidy for $\alpha$, the subgroup $U_0 := \bigcap_{k\in\mathbf{Z}} \alpha^k(U) = U_+\cap U_-$, which is stable under $\alpha$, will appear frequently. Another fact which will be used is that
	\begin{equation}
	\label{eq:tidiness properties}
	\bigcap_{k=0}^n \alpha^k(U) = U\cap \alpha^n(U) = U_+\alpha^n(U) =  \alpha^n(U)U_+.
	\end{equation}

The following is essentially \cite[Lemma 9]{Structure94}, although part \eqref{lem:tidy_criterion2} obtains the same conclusion under a weaker hypothesis and a proof is given for completeness.
\begin{lemma}
	\label{lem:tidy_criterion}
	Let $\alpha\in \Aut(G)$ and suppose that $U$ is tidy for $\alpha$. 
\begin{enumerate}
	\item If $u\in U$ and $\{\alpha^n( u )\}_{n\in\mathbf{N}}$ has an accumulation point, then $u\in U_-$. 	\label{lem:tidy_criterion1}
	\item If $u\in \alpha^{j_1}(U)\dots \alpha^{j_l}(U)$ for some sequence of integers $j_1 < j_2 < \dots< j_l$ and $\{\alpha^n (u)\}_{n\in\mathbf{Z}}$ is contained in a compact set, then $u\in U_0$. 	\label{lem:tidy_criterion2}
\end{enumerate} 
\end{lemma}
\begin{proof}
\eqref{lem:tidy_criterion1} Since $U$ is tidy, $u = u_+u_-$ with $u_\pm\in U_\pm$. If $u_+$ is not in $U_0$, then $\{x^nu_+x^{-n}\}_{n\in\mathbf{N}}$ has no accumulation point. On the other hand, $\{x^nu_-x^{-n}\}_{n\in\mathbf{N}}$ is contained in the compact set $U_-$. Hence, if $\{x^n u x^{-n}\}_{n\in\mathbf{N}}$ has an accumulation point, then $u_+\in U_0$ and $u\in U_-$.

\eqref{lem:tidy_criterion2} Since $U$ is tidy and $j_1 < j_2 < \dots< j_l$, we have that
$$
\alpha^{j_1}(U)\dots \alpha^{j_l}(U)  = \alpha^{j_1}(U_-) \alpha^{j_l}(U_+).
$$
Hence $u = u_-u_+$ with $u_-\in \alpha^{j_1}(U_-)$ and $u_+\in \alpha^{j_l}(U_+)$. Then the argument used in \eqref{lem:tidy_criterion1} shows that $\{\alpha^{-n}(u)\}_{n\in\mathbf{N}}$ is not contained in a compact set if $u_-\not\in U_0$, and $\{\alpha^{n}(u)\}_{n\in\mathbf{N}}$ is not contained in a compact set if $u_+\not\in U_0$. Hence, if $\{x^{n}ux^{-n}\}_{n\in\mathbf{Z}}$ is contained in a compact set, then $u$ belongs to $U_0$.
	\end{proof}


\subsection{Special subgroups}

Minimising, or tidy, subgroups for an automorphism $\alpha$ are not unique in general. Certain other subgroups of $G$, which are uniquely defined in terms of the dynamics of the action of $\alpha$ on $G$, are related to the scale and tidy subgroups. 

Let $\alpha\in\Aut(G)$. The \emph{contraction group for $\alpha$} is 
\begin{align*}
\con{\alpha} &:= \left\{ x\in G \mid \alpha^n(x)\to\ident \text{ as }n\to\infty\right\},\\
\intertext{the \emph{parabolic subgroup} is}
\prb{\alpha} &:= \left\{ x\in G \mid \{\alpha^n(x)\}_{n\in\mathbf{N}} \text{ has compact closure}\right\}\\
\intertext{ and the \emph{Levi subgroup} is}
\lev{\alpha} &:= \prb{\alpha}\cap \prb{\alpha^{-1}} = \left\{ x\in G \mid \{\alpha^n(x)\}_{n\in\mathbf{Z}} \text{ has compact closure}\right\}.\\
\intertext{See \cite[\S3]{ContractionB} for these definitions and the motivations of these names for the subgroups. Contraction subgroups are linked to the scale through \cite[Proposition~3.21]{ContractionB}, which shows that $s(\alpha)$ is equal to the scale of the restriction of $\alpha$ to $\con{\alpha^{-1}}^-$, the closure of the contraction subgroup for $\alpha^{-1}$. Contraction subgroups are not closed in general, but \cite[Theorem~3.32]{ContractionB} shows that $\con{\alpha}$ is closed if and only if the \emph{nub subgroup, $\nub{\alpha}$,} is trivial, where}
\nub{\alpha} &:= \left\{x\in G\mid x\in \con{\alpha}\cap \prb{\alpha^{-1}}\right\}^-.
\end{align*}
In contrast, it is shown in \cite[Proposition~3]{Structure94} that $\prb{\alpha}$, and hence $\lev{\alpha}$ too, is always a closed subgroup of $G$.

Several results about the contraction and Levi subgroups that are used below will now be recalled. The first is easily verified.
\begin{lemma}
	\label{lem:con_alpha^n}
Let $\alpha\in\Aut(G)$. Then $\con{\alpha^n} = \con{\alpha}$ and $\lev{\alpha^n} = \lev{\alpha}$.
\end{lemma}

The next is a corollary of Lemma~\ref{lem:tidy_criterion}.
\begin{lemma}
	\label{lem:levi_stable}
Let $V$ be tidy for $\alpha$ and put $U := V\cap \lev{\alpha}$. Then $\alpha(U) = U$. 
\end{lemma}

The nub subgroup is related to tidiness in \cite[Corollary~4.2]{nub_2014}, which shows that $U$ is tidy below for $\alpha$ if and only if $\nub{\alpha}\leq U$. Combined with  Lemma~\ref{lem:tidy_above} and that fact that $\nub{\alpha}$ is stable under $\alpha$, this implies the following.
\begin{lemma}
	\label{lem:TA+nub-->tidy}
If $U\geq \nub{\alpha}$, then there is $n$ such that $\bigcap_{k=0}^n \alpha^k(U)$ is tidy. 
\end{lemma}
The nub of $\alpha$ is also the largest closed subgroup of $G$ on which $\alpha$ acts ergodically, and the intersection of all subgroups tidy for $\alpha$, as shown in \cite[Theorem~4.1]{nub_2014}. 
 
The proof of the final lemma is a straight-forward verification.
\begin{lemma}
	\label{lem:incl_nub(a)}
Let $\alpha\in\Aut(G)$. Then $\nub{\alpha}\triangleleft \lev{\alpha}$. 
\end{lemma}

The final two lemmas in this subsection use the nub and contraction groups to prove technical facts used in \S\ref{sec:symmetric_relation}.
\begin{lemma}
	\label{lem:include_stableV}
	Let $\alpha\in\Aut(G)$ and $\widetilde{V}$ be a compact $\alpha$-stable subgroup of $G$. Then there is a compact, open subgroup, $V$, of $G$ that is tidy for $\alpha$ and such that $\widetilde{V}\leq V$.
\end{lemma}
\begin{proof}
	The subgroup $\widetilde{V}$ is contained in $\lev{\alpha}$ and so $\widetilde{V}$ normalises $\nub{\alpha}$ by Lemma~\ref{lem:incl_nub(a)}. Hence $\widetilde{V}\nub{\alpha}$ is a compact, $\alpha$-stable subgroup of $G$. Choose a compact, open $U\leq G$ that contains $\widetilde{V}\nub{\alpha}$. (To see that such $U$ exists, choose any compact open subgroup, $W$ say, of $G$ and let $U$ be the product of $\widetilde{V}\nub{\alpha}$ and $\bigcap \{vWv^{-1} \mid v\in \widetilde{V}\nub{\alpha}\}$.) Let ${V} = \bigcap_{k=0}^n \alpha^k(U)$ be the subgroup of $U$ tidy for $\alpha$ shown to exist in Lemma~\ref{lem:TA+nub-->tidy}. Then $\widetilde{V}\leq V$ because $\widetilde{V}$ is $\alpha$-stable.
\end{proof}

\begin{lemma}
	\label{lem:unbounded_sequence}
	Let $\alpha\in\Aut(G)$ and let $v\in \con{\alpha}\setminus \nub{\alpha}$. Then for every $m>0$ and every compact and $\alpha$-stable subgroup $\widetilde{V}$ of $G$, and every sequence $\{w_k\}_{k\in\mathbf{N}}\subset \widetilde{V}$, the sequence
	$$
	\left\{ v\alpha(v)\dots \alpha^k(v)w_k\alpha^{m+k}(v)^{-1} \dots \alpha^{m+1}(v)^{-1}\alpha^m(v)^{-1}\right\}_{k\in\mathbf{N}}
	$$
	is discrete and non-compact.
\end{lemma}
\begin{proof}
Choose $V$ tidy for $\alpha$ and suppose, using Lemma~\ref{lem:include_stableV}, that $\widetilde{V}\leq V$. Replacing $V$ by $\alpha^n(V)$ for some $n\in\mathbf{Z}$, it may be further supposed that $v\in V_+\setminus \alpha^{-1}(V_+)$. Then 
\begin{align*}
& v\alpha(v)\dots \alpha^k(v)w_k\in \alpha^k(V_+)\\
\intertext{and }
& \alpha^{m+k}(v)^{-1} \dots \alpha^{m+1}(v)^{-1}\alpha^m(v)^{-1}\in \alpha^{m+k}(V_+)\setminus \alpha^{m+k-1}(V_+).
\end{align*}
Hence the product is in $\alpha^{m+k}(V)\setminus \alpha^{m+k-1}(V)$ for every $k$ and the claim follows.
\end{proof}

Later sections concern the inner automorphism $\alpha_x : y\mapsto xyx^{-1}$ and the contraction, Levi and nub subgroups are denoted by $\con{x}$, $\lev{x}$ and $\nub{x}$ respectively. 
\subsection{Flat groups and roots}

The main results in this article are about flat groups of automorphisms, that is, groups of automorphisms, $H$, for which there is a compact, open subgroup $U\leq G$ that is minimising for every $\alpha\in H$. As seen in \S\ref{sec:tidy}, this is equivalent to $U$ being tidy for every $\alpha\in H$. 

Many definitions and results about flat groups carry over directly from corresponding ideas for single automorphisms. Of relevance to this article is the subgroup $U_{H0} := \bigcap_{\alpha\in H} \alpha(U)$, defined for any compact open subgroup $U$ but of particular interest if $U$ is tidy. Also, the \emph{lower nub of $H$} is
$$
\lnub{H}:= \left\langle\nub{\alpha} \mid \alpha\in H\right\rangle^-.
$$
The lower nub is defined in \cite{Reid_DynamicsNYJ_2016}, where it is observed that it is contained in every subgroup tidy for $H$, and is therefore compact, but that $\lnub{H}$ may be properly contained in the intersection of the subgroups tidy for $H$. The latter intersection is defined in \cite{Reid_DynamicsNYJ_2016} to be the \emph{nub of $H$}, whence the name \underline{lower} nub for the closed subgroup generated by the nubs of elements of $H$.

The \emph{Levi subgroup for $H$} may also may be defined as follows
$$
\lev{H} := \left\{x\in G \mid \{\alpha(x)\}^-_{\alpha\in H}\text{ is compact}\right\}.
$$
Then $\lev{H}$ is a closed subgroup of $G$, see \cite[Lemma~2.1.8]{flat_Bernoulli}. It is clear that $\lev{H}$ is contained the intersection of all $\lev{\alpha}$ with $\alpha\in H$, but these groups need not be equal, see \cite[Remark~2.1.9]{flat_Bernoulli}. 

Compactness of $\lnub{H}$ and its stability under $H$ imply that $\lnub{H}\leq \lev{H}$. The next lemma then follows from Lemma~\ref{lem:incl_nub(a)}.
\begin{lemma}
	\label{lem:incl_nub}
	Suppose that $H\leq \Aut(G)$ is flat. Then $\lnub{H}\triangleleft \lev{H}$. Hence, if $U$ is a compact subgroup of $\lev{H}$, then $U\lnub{H}$ is a compact subgroup of $\lev{H}$. 
\end{lemma}

Since $U_{H0}$ is compact and $H$-stable, $U_{H0}\leq U\cap \lev{H}$, and  Lemma~\ref{lem:tidy_criterion} implies the reverse inclusion in the following
\begin{lemma}
	\label{lem:flat_tidy_criterion}
	Let $H\leq \Aut(G)$ be flat and suppose that $U$ is tidy for $H$. Then $U\cap \lev{H} = U_{H0}$.
\end{lemma}

It is convenient to state the next couple of results in terms of inner automorphisms and to identify the flat group $H$ with the subgroup of $G$ inducing the inner automorphisms. This is consistent with how they are used later.  The next result is \cite[Theorem~5.9]{SimulTriang} in the case when $K=\triv$ and may be proved in the general case by the argument given for \cite[Proposition~20]{BaumRemyWill}.
\begin{theorem}
	\label{thm:is_flat}
	Suppose that $x_1$, \dots, $x_n\in G$ and the compact group $K\leq G$ satisfy
	\begin{itemize}
		\item $x_iKx_i^{-1} = K$ and
		\item $[x_i,x_j]\in K$.
	\end{itemize}
	Then $\langle x_i, K \mid i\in\{1,\dots,n\}\rangle$ is flat.
\end{theorem}

\begin{definition}
	\label{defn:uniscalar}
	The \emph{uniscalar subgroup}\footnote{The term uniscalar was introduced by T. W. Palmer in \cite[Definition 12.3.25]{PalmerII}.} of the flat group $H$ is
	$$
	H_u := \left\{\alpha\in H \mid s(\alpha)=1=s(\alpha^{-1})\right\} = \left\{\alpha \in H \mid \alpha(V) =V\right\}
	$$
	with $V$ any subgroup minimising for $H$. 
\end{definition}
Then $H_u$ is normal in $H$ and $H/H_u$ is a free abelian group, by \cite[Theorem~6.18]{SimulTriang}. For an alternative proof, see \cite[Theorem~A]{flat_Bernoulli}. The \emph{flat rank} of $H$ is the rank of this free abelian group. The next result is new but, as we shall see, its proof is straightforward.
\begin{proposition}
	\label{prop:frk_equal}
	Let $H = \langle y_1,\dots, y_n\rangle \leq G$ be flat and $V$ be tidy for $H$. Let $L = \langle z_1,\dots, z_n\rangle \leq G$ be another subgroup of $G$ and suppose that for every word $\mathbf{w}(a_1,\dots,a_n)$ in the free group $\langle a_1,\dots, a_n\rangle$, we have 
\begin{equation}
\label{eq:identity}
\mathbf{w}(z_1,\dots,z_n)V	\mathbf{w}(z_1,\dots,z_n)^{-1} = \mathbf{w}(y_1,\dots,y_n)V\mathbf{w}(y_1,\dots,y_n)^{-1}.
\end{equation}
	Then $L$ is flat and $\frk{L} = \frk{H}$.
\end{proposition}
\begin{proof}
	To show that $L$ is flat, it suffices to show that $[zVz^{-1} : zVz^{-1}\cap V] = s(z)$ for every $z\in L$, because that implies that $V$ is minimising for $L$. To this end, let $z\in L$ and suppose that $z = \mathbf{w}(z_1,\dots,z_n)$ with $\mathbf{w}\in \langle a_1,\dots, a_n\rangle$. Put $y = \mathbf{w}(y_1,\dots,y_n)$ in $H$. Then for each $n\in\mathbf{N}$ we have, by Equation~\eqref{eq:identity},
	$$
	[z^nVz^{-n} : z^nVz^{-n}\cap V] = [y^nVy^{-n} : y^nVy^{-n}\cap V],
	$$
and the latter is equal to $s(y)^n$ because $V$ is tidy for $H$. Hence
$$
\lim_{n\to\infty} [z^nVz^{-n} : z^nVz^{-n}\cap V]^{\frac{1}{n}} = s(y)
$$
and it follows, by \cite[Theorem~7.7]{Moller}, that $s(z) = s(y)$. Then, again by Equation~\eqref{eq:identity},
$$
 [zVz^{-1} : zVz^{-1}\cap V] = s(z)
 $$
 and $V$ is minimising for $z$.
 
 The map $\phi : y\mapsto yVy^{-1}$, which sends $H$ to the set of compact open subgroups of $G$, is a bijection $H/H_u\to \phi(H)$. Since, as we now know, $V$ is tidy for $L$ too, the map $\psi : z\mapsto zVz^{-1}$ is also a bijection $L/L_u \to \psi(L)$. Equation~\eqref{eq:identity} implies that $\phi(H) = \psi(L)$ and so $\phi$ and $\psi$ induce a bijection $H/H_u\to L/L_u$. Moreover, if $\mathbf{w}_1$ and $\mathbf{w}_2$ are two words in the free group $\langle a_1,\dots, a_n\rangle$, and $y_i = \mathbf{w}_i(y_1,\dots,y_n)$ and $z_i = \mathbf{w}_i(z_1,\dots, z_n)$ are the corresponding elements of $H$ and $L$, then~\eqref{eq:identity} implies that $\phi(y_1y_2) = \psi(z_1z_2)$ and hence that the bijection $H/H_u \to L/L_u$ is a group isomorphism. The free abelian groups $H/H_u$ and $L/L_u$ therefore have the same rank.
\end{proof}

\section{The Levi subgroup and tidy subgroups for a flat group}
\label{sec:Levi}

Preliminary results on links between the Levi subgroup and subgroups tidy for a flat group are derived in this section. They will be used in the next section in the proof of the main theorem. Many of these results may have independent interest. The first shows that the property of being tidy for a flat group $H$ is virtually invariant under conjugation by elements of $\lev{H}$. 
\begin{proposition}
	\label{prop:Htilde}
	Let $H\leq G$ be flat, $y\in \lev{H}$ and $V$ be a compact open subgroup of $G$ tidy for $H$. Then $$
L := \left\{x\in H\mid xyx^{-1}\in yV \right\} 
	$$
 is a finite index subgroup of $H$ and $yVy^{-1}$ is tidy for $L$.
\end{proposition}
\begin{proof}
Suppose that $x\in L$. Then $xyx^{-1} = yv$ with $v\in V$ and, for each $h\in H$,
\begin{equation}
\label{eq:Htilde}
h(xyx^{-1})h^{-1} = h(yv)h^{-1} = (hyh^{-1}) (hvh^{-1})
\end{equation}
for all $n\in\mathbf{Z}$. Hence $hvh^{-1}$  is in the compact set $(y^H)^2$ for all $n\in\mathbf{Z}$ and, since $V$ is tidy for $H$, Lemma~\ref{lem:flat_tidy_criterion} implies that  $v\in V_{H0}$. Let $x_1,x_2\in L$. Then $x_1yx_1^{-1} = yv_1$ and $x_2^{-1}yx_2 = yv_2$ with $v_1,v_2\in V_{H0}$ and  
$$
(x_1x_2^{-1})y(x_1x_2^{-1})^{-1} = x_1 yv_2 x_1^{-1} = y(v_1 x_1 v_2x_1^{-1}) 
$$
with $v_1 x_1 v_2x_1^{-1} \in V_{H0}$. Hence $x_1x_2^{-1}\in L$ and $L$ is a group.

Since $y\in\lev{H}$, $\{xyx^{-1}V\mid x\in H\}$ is a finite set of $V$-cosets and it follows, again using Equation~\eqref{eq:Htilde}, that $H$ is covered by a finite number of $L$ cosets. 

For each $x\in L$, there is $v\in V_{H0}$ such that
$$
x(yVy^{-1})x^{-1} = (xyx^{-1})(xVx^{-1})(xyx^{-1})^{-1} = yv(xVx^{-1})v^{-1}y^{-1} = y(xVx^{-1})y^{-1}.
$$
Hence
\begin{align*}
&\phantom{=}\ \,[x(yVy^{-1})x^{-1} : x(yVy^{-1})x^{-1}\cap yVy^{-1}] \\
&= [y(xVx^{-1})y^{-1}:y(xVx^{-1})y^{-1}\cap yVy^{-1}]\\
&= [xVx^{-1}:xVx^{-1}\cap V] = s(x)
\end{align*}
and $yVy^{-1}$ is tidy for $x$. Therefore $yVy^{-1}$ is tidy for $L$. 
\end{proof}

The proof of Proposition~\ref{prop:Htilde} in fact shows the following stonger assertion.
\begin{corollary}
	\label{cor:Htilde}
	Let $y$ and $L$ be as in Proposition~\ref{prop:Htilde}. Then $xyx^{-1}y^{-1}\in V_{H0}$ for every $x\in L$.
\end{corollary}

In the next section, a flat group $H$ will be augmented by adding an element from its Levi subgroup. The following proposition helps to find a subgroup tidy this element.
\begin{proposition}
	\label{prop:nub_in_lev}
	Let $H\leq G$ be flat and $y\in\lev{H}$. Then $\nub{y}\leq \lev{H}$.
	\end{proposition}
\begin{proof}
Choose $V$ be tidy for $H$. Then Lemma~\ref{lem:tidy_above} shows that there is $n\geq0$ such that $V_1 := \bigcap_{k=0}^n y^kVy^{-k}$ is tidy above for $y$; and Proposition~\ref{prop:Htilde} that there is a finite index subgroup, $H_1$, of $H$ such that $V_1$ is tidy for $H_1$. Set $V_2 = y^{-1}V_1y\cap V_1$. Then, by Proposition~\ref{prop:Htilde} again, there is a finite index subgroup, $H_2$, of $H_1$ such that $V_2$ is tidy for $H_2$ and $H_3 := \left\{x\in H_2\mid xyx^{-1}\in yV _2\right\}$ has finite index in $H_2$. Consider $c\in\nub{y}$ and let $x\in H_3$. Then $xcx^{-1}\in \nub{xyx^{-1}}$. We have $xyx^{-1}\in yV_2$ because $x\in H_3$ and so, by \cite[Corollary 3.4]{CapReidW_Titscore}, there is $r\in V_1$ such that $\nub{xyx^{-1}} = r\nub{y}r^{-1}$. Hence $xcx^{-1}$ is in the compact set $V_1\nub{y}V_1$ for every $x\in H_3$. Since $H_3$ has finite index in $H$, by Proposition~\ref{prop:Htilde}, it follows that $\{xcx^{-1}\}_{x\in H}$ has compact closure and hence that $c\in \lev{H}$.
\end{proof}
The previous argument shows that, if it is supposed only that $c\in\con{y}$ (rather than $\nub{y}$) and $x\in H_3$, then $xcx^{-1}$ belongs to $t\con{y}t^{-1}$ for some $t\in V$, by \cite[Corollary 3.2]{CapReidW_Titscore}. However, the remainder of the argument does not apply and it is not true in general  that $\con{y}\leq\lev{H}$, as may be seen by considering $H = \langle y\rangle$ for example.

Augmenting a flat group $H$ involves modifying a subgroup, $U$, tidy for $H$. The following lemma facilitates one such modification. 
\begin{lemma}
	\label{lem:nub_in_tidy}
Suppose that $U$ is tidy for $x\in G$ and that $C$ is a compact subgroup of $G$ such that $CU = UC$ and $\{x^nCx^{-n}\}_{n\in\mathbf{Z}}\subseteq (xUx^{-1})C$. Then $UC$ is a compact open subgroup of $G$ that is tidy for $x$. 
\end{lemma}
\begin{proof}
	Since $U$ is compact and open and $C$ is compact, $U$, $UC$ is a compact open subset of $G$ and the condition $UC=CU$ implies that it is a group. Tidiness of $UC$ for $x$ follows from the fact that 
	$$
	[x(UC)x^{-1} : x(UC)x^{-1}\cap UC] = [xUx^{-1} : xUx^{-1}\cap U] = s(x).
	$$
	To see this, consider the map 
	$$
	xUx^{-1}/(xUx^{-1}\cap U) \to x(UC)x^{-1} / (x(UC)x^{-1}\cap UC)
	$$ 
	defined by $w(xUx^{-1}\cap U) \mapsto w(x(UC)x^{-1}\cap UC)$. This map is well-defined because $xUx^{-1}\cap U\leq x(UC)x^{-1}\cap UC$ and is onto because 
	$$
	x(UC)x^{-1} = (xUx^{-1})(xCx^{-1}) \subseteq (xUx^{-1})C.
	$$ 
	The map is one-to-one because, if 
	$$
	w_1,w_2\in xUx^{-1}\text{ and }w_2^{-1}w_1 \in x(UC)x^{-1}\cap UC,
	$$ 
	then $w_2^{-1}w_1 = wc$ with $w\in U$ and $c\in C$. Hence $c = w_2^{-1}w_2w^{-1}$, which belongs to $(xUx^{-1})U$, and $\{x^ncx^{-n}\}_{n\in \mathbf{Z}}$ has compact closure because it is a subset of $xUx^{-1}C$. Then, since $U$ is tidy for $x$, Lemma~\ref{lem:tidy_criterion}\eqref{lem:tidy_criterion2} shows that $w_2^{-1}w_1w^{-1} =: w_0$ belongs to $U_0$, and hence that $w_2^{-1}w_1 = w_0w\in U\cap xUx^{-1}$. 
\end{proof}

\begin{corollary}
	\label{cor:nub_in_tidy}
	The	hypotheses of Lemma~\ref{lem:nub_in_tidy} imply that $\{y^nCy^{-n}\}_{n\in\mathbf{Z}}\subseteq U_0C$.
\end{corollary}
\begin{proof}
	The	hypotheses imply that, if $c\in C$ and $k\in\mathbf{Z}$, then $y^kcy^{-k} = uc_1$ with $u\in yUy^{-1}$ and $c_1\in C$. Then 
	$$
	\{y^nuy^{-n}\}_{n\in\mathbf{Z}} = \{y^n(y^kcy^{-k}c_1^{-1})y^{-n}\}_{n\in\mathbf{Z}} \subseteq ((yUy^{-1})C)^2,
	$$ 
	which is compact. Hence $u\in U_0$ by  Lemma~\ref{lem:tidy_criterion}\eqref{lem:tidy_criterion2}. 
\end{proof}

The following proposition shows that a subgroup tidy for $y\in\lev{H}$ and $H$ may be found at the expense passing to a finite index subgroup, $F$, of $H$. Note that this proposition does not imply that the subgroup found is tidy for $\langle y, F\rangle$ however. That is done in the next section under additional hypotheses. 
\begin{proposition}
	\label{prop:tidy_for_both}
	Let $H\leq G$ be flat and $y\in \lev{H}$. Then there are a compact open subgroup, $V$, of $G$ and a finite index subgroup, $F\leq H$, such that $V$ is tidy for $F$ as well as for $y$.
\end{proposition}
\begin{proof}
Let $W$ be a subgroup of $G$ tidy for $H$. Then, since $\nub{y}$ is compact and $W$ is open, there are $c_i\in\nub{y}$, $i\in\{1,\dots,k\}$, such that 
$$
\bigcap\left\{cWc^{-1}\mid c\in\nub{y}\right\} = \bigcap_{i=1}^k c_iWc_i^{-1}.
$$ 
By Lemma~\ref{prop:nub_in_lev}, each $c_i\in\lev{H}$ and hence, by Proposition~\ref{prop:Htilde}, there is a finite index subgroup $H_i\leq H$ such that $c_iWc_i^{-1}$ is tidy for $H_i$. Then $W_1 := \bigcap_{i=1}^k c_iWc_i^{-1}$ is tidy for $\bigcap_{i=1}^k H_i$ and is normalised by $\nub{y}$. Set $W_2 = \nub{y}W_1$. Then $W_2$ is a compact open subgroup of $G$ that contains $\nub{y}$.

We shall see that $W_2$ is tidy for a finite index subgroup of $\bigcap_{i=1}^k H_i$. Since $\nub{y}\leq \lev{H}$, by Proposition~\ref{prop:nub_in_lev}, and is compact, there is, by Proposition~\ref{prop:Htilde}, a finite index subgroup, $E$, of $\bigcap_{i=1}^k H_i$ such that $x\nub{y}x^{-1}\subseteq \nub{y}W_1$ for every $x\in E$. Corollary~\ref{cor:Htilde} shows that, in fact, $x\nub{y}x^{-1}\subseteq \nub{y}(W_1)_{E0}$ for every $x\in E$. Therefore, taking $U = W_1$, $C = \nub{y}$, the hypotheses of Lemma~\ref{lem:nub_in_tidy} are satisfied for every $x\in E$ and we conclude that $W_2$ is tidy for every such $x$. 

Finally, choose $n$ sufficiently large that  $V := \bigcap_{k=0}^n y^kW_2y^{-k}$ is tidy above for~$y$. Such $n$ exists by Lemma~\ref{lem:tidy_above}. Then, since $\nub{y}\leq W_2$, Lemma~\ref{lem:TA+nub-->tidy} shows that $V$ is tidy for $y$. By Proposition~\ref{prop:Htilde}, there is for each $k\in\{0,\dots,n\}$ a finite index subgroup, $E_k$, of $E$ such that $y^kW_2y^{-k}$ is tidy for $E_k$. Then $V$ is tidy for the subgroup $F:= \bigcap_{k=0}^n E_k$ which has finite index in $E$, and hence in $H$. 
\end{proof}

The final proposition in this section derives a useful property of subgroups tidy above for an element of $G$. 
\begin{proposition}
	\label{prop:conj_into_U0}
	Let $U$ be tidy above for $y$ and suppose that $z\in UyU$. Then there is $u\in U$ such that $uzu^{-1}\in yU_0$.
\end{proposition}
\begin{proof}
	Let $v,w\in U$ be such that $z = vyw$. Since $U$ is tidy above for $y$, we have that $v=v_-v_+$ with $v_\pm\in U_\pm$. Then $vyw = v_-yw'$ with $w' = y^{-1}v_+yw\in U$. Repeating the argument for $w'$ shows that it may be supposed that 
	$$
	z = vyw\text{ with }v\in U_-\text{ and }w\in U_+.
	$$
	
Next, setting $v_0=v$ and $w_0=w$, we recursively construct $v_i\in y^{i}U_-y^{-i}$ and $w_i\in U_+$ such that 
	$$
	v_i^{-1}(v_iyw_i)v_i = v_{i+1}yw_{i+1}.
	$$
	Suppose that $v_i$ and $w_i$ have been constructed for some $i$. Then
	$$
		v_i^{-1}(v_iyw_i)v_i = yw_iv_i = yv_i'w_i'
	$$
	with $v_i'\in y^{i}U_-y^{-i}$ and $w_i\in U_+$ because $w_iv_i \in U_+ \left(y^iU_-y^{-i}\right) = \left(y^iU_-y^{-i}\right)U_+ $ by Equation~\eqref{eq:tidiness properties}. Setting $v_{i+1} = yv_i'y^{-1}$ and $w_{i+1}=w_i'$ continues the construction. Let $u_i = v_0v_1\dots v_{i-1}$. Then for each $i\geq1$ we have
	$$
	u_i(vyw)u_i^{-1} = v_iyw_i
	$$
with $v_i\in y^kU_-y^{-k}$ and $w_i\in U_+$. Since $U$, $U_-$ and $U_+$ are compact, the sequences $\{u_i\}_{i\in\mathbf{Z}}$, $\{v_i\}_{i\in\mathbf{Z}}$ and $\{w_i\}_{i\in\mathbf{Z}}$ have limit points $\hat{u}\in U$, $\hat{v}\in \bigcap_{k\geq0} y^kU_-y^{-k} = U_0$ and $\hat{w}\in U_+$ respectively with $\hat{u}(vyw)\hat{u}^{-1} = \hat{v}y\hat{w}$. This may be further simplified to
$$
\hat{u}(vyw)\hat{u}^{-1} = y\hat{w}
$$
with $\hat{w}\in U_+$ because $y^{-1}\hat{v}y \in U_0$ and $U_0\leq U_+$.

Repeating the construction of the previous paragraph: starting with $\hat{w}_0 = \hat{w}$ and defining $\hat{w}_{i+1} = y^{-1}\hat{w}_iy$ for $i\geq0$; and then defining $\hat{u}_i = \hat{w}_{i-1}\dots \hat{w}_0$ for $i\geq1$ yields
$$
\hat{u}_i\hat{u}(vyw)\hat{u}^{-1}\hat{u}_i^{-1} = y\hat{w}_i
$$
with $\hat{w}_i\in y^{-i}U_+y^i$. Since $U_+$ is compact, the sequences $\{\hat{u}_i\}_{i\geq1}$ and $\{\hat{w}_i\}_{i\geq1}$ have limit points $\doublehat{u}\in U$ and $u_0\in \bigcap_{i\geq0} y^{-i}U_+y^i = U_0$ respectively such that
$$
\doublehat{u}\hat{u}(vyw)\hat{u}^{-1}\doublehat{u}^{-1} = yu_0.
$$
Therefore the claim holds with $u = \doublehat{u}\hat{u}$. 
\end{proof}

\section{Groups with flat rank greater than $1$}
\label{sec:Main}

In this section it is seen that a flat subgroup, $H$, of $G$ for which there is $y\in\lev{H}$ with scale greater than~$1$ may, under certain hypotheses, be modified to produce a subgroup with flat rank increased by $1$. One example is given that illustrates the method and another that shows its limitations. 

The hyptheses are that $H$ has the form $\langle x_i, K \mid i\in\{1,\dots,n\}\rangle$ described in Theorem~\ref{thm:is_flat}. The theorem shows how to produce a group $\langle y, z_i, L\mid i\in\{1,\dots,n\}\rangle$ that is also flat and may have rank greater than that of $H$. The group produced by the theorem also has the required form, as does the singly generated group $\langle y\rangle$. Hence, the method may be iterated to produce groups with higher flat rank so long as there is a non-uniscalar element in the Levi subgroup of the flat group last produced. In the statement of the theorem, $s_{\lev{H}}(y)$ denotes the scale of $y$ in the group $\lev{H}$, which may be less than the scale of $y$ in $G$. 
\begin{theorem}
	\label{thm:extend_flat}
	Let $G$ be a \tdlc~group and let $x_i\in G$, $i\in\{1,\dots,n\}$, and $K$ be a compact subgroup of $G$ be such that: 
	\begin{itemize}
		\item $x_iKx_i^{-1} = K$; and
		\item $[x_i,x_j]\in K$. 
	\end{itemize}
Denote the flat group $\langle x_i, K \mid i\in\{1,\dots,n\}\rangle$ by $H$ and suppose that $y\in\lev{H}$. 

Then there are $z_i\in G$ and a compact group $L\leq G$ such that for all $i\in\{1,\dots,n\}$: 
\begin{enumerate}
	\item $z_iLz_i^{-1} = L$, $yLy^{-1}=L$; 
\label{thm:extend_flat2}	
	\item $[z_i,z_j]\in L$ and $[z_i,y]\in L$; and
\label{thm:extend_flat3}
	\item $\con{z_i}=\con{x_i}$ and $\lev{z_i}=\lev{x_i}$
\label{thm:extend_flat1}
\end{enumerate}
Hence $\langle y, z_i, L\mid i\in\{1,\dots,n\}\rangle$ is flat. If $s_{\lev{H}}(y)>1$, then 
$$
\frk{\langle y, z_i, L\mid i\in\{1,\dots,n\}\rangle} = \frk{\langle x_i, K \mid i\in\{1,\dots,n\}\rangle}+1.
$$
\end{theorem}	
\begin{proof}	
Proposition~\ref{prop:tidy_for_both} shows that there is $V\leq G$ compact, open and tidy both for $y$ and a finite index subgroup, $F$, of $H$. Let $U := V\cap \lev{H}$. Then $U$ is a compact open subgroup of $\lev{H}$ that is tidy for $y\in\lev{H}$, and Lemma~\ref{lem:levi_stable} shows that $xUx^{-1} = U$ for every $x\in F$. 

We establish \eqref{thm:extend_flat2} and \eqref{thm:extend_flat3} through a series of claims. For these, note that the hypotheses that $[x_i,x_j]\in K$ and that $K$ is a compact subgroup normalised by each $x_i$ imply that $H\leq \lev{H}$. As a consequence of this, $z_i$ and $L$ may be defined within $\lev{H}$. With this in mind, it will be convenient to suppose that $G = \lev{H}$, and hence that $U=V$ and $H$ is uniscalar, while proving \eqref{thm:extend_flat2} and \eqref{thm:extend_flat3}. 
	\begin{Claim}
		\label{claim:exists_mi}
	For each $i\in\{1,\dots,n\}$ there is $n_i>0$ such that 
	$$
	x_i^{n_i}Ux_i^{-n_i} = U \text{ and } x_i^{n_i}yx_i^{-n_i}\in yU.
	$$
	\end{Claim}
\begin{proof}
Since $F$ has finite index in $H$, there is $m_i>0$ such that $x_i^{m_i}\in F$ and hence that $x_i^{m_i}Ux_i^{-m_i} = U$. Proposition~\ref{prop:Htilde} shows that there is a finite index subgroup $L\leq \langle x_i^{m_i} \mid i\in\{1,\dots,n\}\rangle$ such that $xyx^{-1}\in yU$ for every $x\in L$. Hence  there is $n_i$, a multiple of $m_i$, such that $	x_i^{n_i}Ux_i^{-n_i} = U$ and $x_i^{n_i}yx_i^{-n_i}\in yU$. 
	\end{proof}

Next, Proposition~\ref{prop:conj_into_U0} shows that for each $i\in\{1,\dots,n\}$ there is $u_i\in U$ such that 
$$
u_ix_i^{n_i}yx_i^{-n_i}u_i^{-1}\in yU_0.
$$
Put $z_i = u_ix_i^{n_i}$. Then it follows from this and Claim~\ref{claim:exists_mi}.~that
\begin{Claim}
	\label{claim:zy_commutator}
	$z_iyz_i^{-1}\in yU_0$ and $z_iUz_i^{-1} = U$ for every $i\in\{1,\dots,n\}$.
\end{Claim}

Having defined the elements $z_i$, the subgroup $L$ will be defined next. Towards this, let $K'$ be the smallest closed subgroup satisfying, for all $i,j,k$:
\begin{itemize}
	\item $[x_i^{n_i},x_j^{n_j}]\in K'$ and 
	\item $x_k^{n_k}K'x_k^{-n_k} = K'$.
\end{itemize}
Then $K'$ is densely generated by a normal subgroup of $\langle x_i^{n_i}\mid i\in\{1,\dots,n\}\rangle$ and so $wUw^{-1} = U$ for all $w\in K'$ by Claim~\ref{claim:exists_mi}. Furthermore, $K'$ is contained in $K$ and therefore is compact. Hence $UK'$ is a compact group stable under conjugation by all $x_i^{n_i}$. Since $U$ is normalised by all $x_i^{n_i}$, $UK'$ is also stable under conjugation by all $z_k$ and $[z_i,z_j] = [u_ix_i^{n_i},u_jx_j^{n_j}]\in UK'$ for all $i,j,k\in\{1,\dots,n\}$. Hence 
\begin{Claim}
	\label{claim:comm_in_K''}
The smallest closed subgroup, $K''$, satisfying that $[z_i,z_j]\in K''$ and $z_kK''z_k^{-1} = K''$ 
 for all $i,j,k\in\{1,\dots,n\}$ is a compact subgroup of $UK'$. 
\end{Claim}

It follows from Claim~\ref{claim:zy_commutator}.~and invariance of $U_0$ under conjugation by $y$ that $z_iy^lz_i^{-1}\in y^lU_0$ and hence that $z_i(y^lUy^{-l})z_i^{-1} = y^lUy^{-l}$ for every $l\in\mathbf{Z}$. Hence $z_iU_0z_i^{-1} = U_0$ for each $i\in\{1,\dots,n\}$ and  $zU_0z^{-1} = U_0$ for every $z\in K''$. Therefore $K''U_0$ is a group, which is compact because $K''$ and $U_0$ are. Define $L = K''U_0$. We are now ready to prove parts	\eqref{thm:extend_flat2} and \eqref{thm:extend_flat3}.
\begin{Claim}
	\label{claim:L}
For every $i,j\in\{1,\dots,n\}$
	$$
	[z_i,y]\in L,\ [z_i,z_j]\in L\text{ and }yLy^{-1} = L = z_iLz_i^{-1}.
$$  
\end{Claim}
\begin{proof}[Proof of claim] That $[z_i,y]\in L$ is because $[z_i,y]\in U_0$ by Claim~\ref{claim:zy_commutator}., and that $[z_i,z_j]\in L$ is because $[z_i,z_j]\in K''$ by definition, see Claim~\ref{claim:comm_in_K''}. Since, as was just shown, $z_iU_0z_i^{-1} = U_0$ and,  by definition, $z_iK''z_i^{-1} = K''$, it follows that $z_iLz_i^{-1} = L$ for each $i\in\{1,\dots,n\}$. The definition of $U_0$ implies that $yU_0y^{-1} = U_0$. Finally, $K''$ is generated by elements of the form $z_{k_1}\dots z_{k_r}[z_i,z_j]z_{k_r}^{-1}\dots  z_{k_1}^{-1}$ and for any such element we have
	$$
	y(z_{k_1}\dots z_{k_r}[z_i,z_j]z_{k_r}^{-1}\dots  z_{k_1}^{-1})y^{-1}\in z_{k_1}\dots z_{k_r}[z_i,z_j]z_{k_r}^{-1}\dots  z_{k_1}^{-1} U_0
	$$
because $[z_i,y]\in U_0$ by Claim \ref{claim:zy_commutator}.~and $U_0$ is normalised $z_1$, \dots, $z_n$. Therefore $yK''y^{-1}\leq K''U_0 = L$.
	\end{proof}

The remainder of the proof returns to the general case where $G\ne \lev{H}$ and $V\ne U$. To see that $\con{z_i} = \con{x_i}$, recall from Lemma~\ref{lem:con_alpha^n} that $\con{x_i^{n_i}} = \con{x_i}$ and consider $c\in G$. Then,  since $x^{n_i}$ normalises $U$, we have for each $r\in\mathbf{N}$ that
$$
z_i^rcz_i^{-r} = (u_ix_i^{n_i})^rc(u_ix_i^{n_i})^{-r} = t_r (x_i^{rn_i}cx_i^{-rn_i})t_r^{-1}
$$ 
with $t_r\in U$. Since the open normal subgroups of $V$ form a base of neighbourhoods of $\ident$, it follows that $c\in \con{z_i}$ if and only if $c\in\con{x_i^{n_i}}$  . That $\lev{z_i} = \lev{x_i}$ holds by the same argument and the fact that $\lev{x_i^{n_i}} = \lev{x_i}$.
	
The claim that $\langle y, z_i, L \mid i\in\{1,\dots,n\}\rangle$ is flat follows from Theorem~\ref{thm:is_flat} and~\eqref{thm:extend_flat2} and~\eqref{thm:extend_flat3}. The statement about the flat ranks will follow from 
\begin{Claim}
	\label{claim:fltrk}
	$\frk{\langle z_i,L\mid i\in\{1,\dots,n\}\rangle} = \frk{H}$.
\end{Claim}
\begin{proof}  Since $V\leq G$ is tidy for $\langle x_i^{n_i},\mid i\in\{1,\dots,n\}\rangle$, we have, for each word $\mathbf{w}(a_1,\dots,a_n)$ in the free group $\langle a_1,\dots, a_n\rangle$, 
	$$
	\mathbf{w}(z_1,\dots,z_n) = \mathbf{w}(x_1^{n_1},\dots,x_n^{n_n})u
	$$ 
	with $u\in U$ because each $x_i^{n_i}$ normalises $U$ for each $i\in\{1,\dots,n\}$. Hence
	$$
	\mathbf{w}(z_1,\dots,z_n)U	\mathbf{w}(z_1,\dots,z_n)^{-1} = \mathbf{w}(x_1^{n_1},\dots,x_n^{n_n})U\mathbf{w}(x_1^{n_1},\dots,x_n^{n_n})^{-1}.
$$
Therefore, by Proposition~\ref{prop:frk_equal}, 
$$
\frk{\langle z_i\mid i\in\{1,\dots,n\}\rangle} = \frk{\langle x_i^{n_i}\mid i\in\{1,\dots,n\}\rangle}.
$$ 
This suffices to prove the claim because $L$ and $K$ are contained in the uniscalar subgroups of the flat groups $\langle z_i,L\mid i\in\{1,\dots,n\}\rangle$ and $H$.
\end{proof}

Now we are ready to show that the flat rank increases by $1$ if $s_{\lev{H}}(y)>1$. If that is the case, then $U_{y++} := \bigcup_{k\in\mathbf{Z}} y^k U_{y+}y^{-k}$ is a closed, non-compact subgroup of $\lev{H}$ because $U$ is tidy for $y$ in $\lev{H}$. Furthermore, since $y^kUy^{-k}$ is stable under conjugation by $z_i$, $U_{y++}$ is also stable under conjugation by $z_i$ for each $i\in\{1,\dots,n\}$. Therefore $\langle y, z_i\mid i\in\{1,\dots,n\}\rangle$ acts on $U_{y++}$ by conjugation. The map $\alpha_x : u\mapsto xux^{-1}$ is an automorphism of $U_{y++}$, and 
$$
x\mapsto \alpha_x : \langle y, z_i\mid i\in\{1,\dots,n\}\rangle \to \Aut(U_{y++})
$$
is a homomorphism. Denoting the modular function on $\Aut(U_{y++})$ by $\Delta$, define $\rho : \langle y, z_i\mid i\in\{1,\dots,n\}\rangle \to (\mathbf{R},+)$ by 
$$
\rho(x) = \log\Delta(\alpha_x), \qquad x\in \langle y, z_i\mid i\in\{1,\dots,n\}\rangle.
$$ 
Then $\rho$ is a homomorphism with $\langle z_i\mid i\in\{1,\dots,n\}\rangle\leq \ker\rho$ because $U_{y++}$ is contained in $\lev{H}$, which is equal to $\lev{\langle z_i\mid i\in\{1,\dots,n\}\rangle}$. On the other hand, $\rho$ maps onto a subgroup of $(\mathbf{R},+)$ isomorphic to $(\mathbf{Z},+)$ because $y$ is not uniscalar. Therefore the rank of $\langle y,z_i,L\mid i\in\{1,\dots,n\}\rangle$ modulo its uniscalar subgroup is $\frk{H}+1$. 
\end{proof}

Since finitely generated abelian groups are flat, by Theorem~\ref{thm:is_flat}, Theorem~\ref{thm:extend_flat} specialises to the following.
\begin{corollary}
Let $G$	be a \tdlc~group and $x_i\in G$, $i\in\{1,\dots, n\}$ commute and denote $H = \langle x_1,\dots, x_n\rangle$. Suppose that $y\in C(H)$. Then $\langle H,y\rangle$ is flat. If $s_{C(H)}(y)>1$, then the flat-rank of $\langle H,y\rangle$ is equal to the flat-rank of $H$ plus~$1$.
\end{corollary}

The following example illustrates the steps in the proof of Theorem~\ref{thm:extend_flat}.
\begin{example}
	\label{examp:illustrate}
	Let $G = GL(3,\mathbf{Q}_p)$ and $H = \langle x\rangle$ with 
	$$
	x = \left(\begin{matrix}
	1 & a & 0\\ 0 & 1 & 0\\ 0 & 0 & p
	\end{matrix}\right) \text{ and }a\in\mathbf{Q}_p\setminus\{0\}.
	$$
(We will use that $|a|_p$, the $p$-adic absolute value of $a$, is a power of $p$.) Then $H$ is flat and
	$$
	\lev{H} = \left\{\left(\begin{matrix}
	a_{11} & a_{12} & 0\\ a_{21} & a_{22} & 0\\ 0 & 0 & a_{33}
	\end{matrix}\right) \mid a_{ij}\in\mathbf{Q}_p,\ (a_{11}a_{22}-a_{12}a_{21})a_{33}  \ne 0\right\}.
	$$
Since $s(x) = p^2$, $H$ has flat rank equal to $1$.
	
Hence
$$
y = \left(\begin{matrix}
p & 0 & 0\\ 0 & 1 & 0\\ 0 & 0 & 1
\end{matrix}\right)
$$
belongs to $\lev{H}$ and $s_{\lev{H}}(y) = p$. 

We apply the steps of the construction to this $H$ and $y$. To begin, the group
$$
V = \left\{\begin{pmatrix}
a_{11} & pa_{12} & pa_{13}\\ a_{21} & a_{22} & pa_{23}\\ a_{31} & a_{32} & a_{33}
\end{pmatrix}\in GL(3,\mathbf{Q}_p)\mid a_{ij}\in\mathbf{Z}_p,\ |a_{11}a_{22}a_{33}|_p=1\right\}
$$
is tidy for both $H$ and $y$ and 
$$
U = V\cap \lev{H} = \left\{\begin{pmatrix}
a_{11} & pa_{12} & 0\\ a_{21} & a_{22} & 0\\ 0 & 0 & 1
\end{pmatrix}\in GL(3,\mathbf{Q}_p)\mid a_{ij}\in\mathbf{Z}_p,\ |a_{11}a_{22}|_p=1\right\}.
$$
Then choosing $n = p^2|a^{-1}|_p$ ensures that $x^nUx^{-n} = U$ and $x^nyx^{-n}\in yU$ as in Claim~\ref{claim:exists_mi}. Next, put 
$$
u = \begin{pmatrix}
1 & -na & 0\\ 0 & 1 & 0\\ 0 & 0 & 1
\end{pmatrix}.
$$
Then $u\in U$ and $ux^nyx^{-n}u^{-1} = y$. Put  $z = ux^n$. Then
$$
z = \begin{pmatrix}
1 & 0 & 0\\ 0 & 1 & 0\\ 0 & 0 & p^n
\end{pmatrix},
$$
and Claim~\ref{claim:zy_commutator} holds and the subgroups $K''$ and $L$ in Claims~\ref{claim:comm_in_K''} and~\ref{claim:L} are trivial because $y$ and $z$ commute. Hence $\langle y, z\rangle$ is flat and has flat rank equal to $2$. 
\end{example}

The proof of Theorem~\ref{thm:extend_flat} and construction in Example~\ref{examp:illustrate} pass from the given flat group $H$ to a finite index subgroup on the way to constructing the new flat group with increased rank. That could have been avoided in Example~\ref{examp:illustrate} by choosing a different subgroup $V$ tidy for $H$ and $y$. The next example shows that passing to a finite index subgroup cannot always be avoided.
\begin{example}
	Define $x,y\in \mathbf{Z}^2\to \Aut(SL(2,\mathbf{Q}_p)\times \mathbf{Q}_p)$  by 
	\begin{align*}
	x(A,c) &= \left(\begin{pmatrix}
	0&1\\1&0
	\end{pmatrix}A\begin{pmatrix}
	0&1\\1&0
	\end{pmatrix},pc\right)\\
\text{ and }	y(A,c) &= \left(\begin{pmatrix}
	p&0\\0&1
	\end{pmatrix}A\begin{pmatrix}
	p&0\\0&1
	\end{pmatrix}^{-1}, c\right),\qquad A\in SL(2,\mathbf{Q}_p),\ c\in \mathbf{Q}_p,
	\end{align*}
	and let $G = \langle x,y\rangle\ltimes \left(SL(2,\mathbf{Q}_p)\times \mathbf{Q}_p)\right)$. Then $\langle x^2,y\rangle$ is finitely generated and abelian, and hence flat. In fact, if $U\leq SL(2,\mathbf{Q}_p)$ is the subgroup 
	$$
	\left\{\begin{pmatrix} a& pb\\ c&d \end{pmatrix}\mid a,b,c,d\in \mathbf{Z}_p,\ ad-pbc =1\right\},
	$$
	then $U\times \mathbf{Z}_p$ is tidy for $\langle x^2,y\rangle$. However, $\langle x,y\rangle$ is not flat because, for example, its abelianisation is finite.
\end{example}

The next example shows that the construction in Theorem~\ref{thm:extend_flat} might not apply because it may happen that, while $\lev{H}$ contains elements of a flat group with larger flat rank, it does not contain a contraction subgroup for any element of the larger flat group. 
\begin{example}
	Let $x_1$ and $x_2$ be the automorphisms of $\mathbf{Q}_p^2$ given by
	$$
	x_1(\xi_1,\xi_2) = (p\xi_1,p\xi_2)\text{ and }x_2(\xi_1,\xi_2) = (p\xi_1,p^{-1}\xi_2),\qquad (\xi_1,\xi_2)\in\mathbf{Q}_p^2,
	$$ 
	and let $A := \langle x_1,x_2\rangle$. Then $A$ is a flat subgroup of $G := \mathbf{Q}_p^2\rtimes A$ having $V := \mathbf{Z}_p^2\rtimes \{0\}$ as a tidy subgroup. The uniscalar subgroup of $A$ is trivial and so the flat rank of $A$ is equal to the rank of $A$, which is $2$.
	
	Let $H := \langle x_1\rangle \leq A$. Then the flat rank of $H$ is equal to $1$. However, $\lev{H} = A$ which is discrete and hence uniscalar, that is, $s_{\lev{H}}(y)=1$ for every $y\in A$. Hence no $y$ can be chosen from $\lev{H}$ that can be used in the construction in Theorem~\ref{thm:extend_flat}.
\end{example}

The final example shows that the construction in Theorem~\ref{thm:extend_flat} might not apply because it may happen that, while $\lev{H}$ contains a contraction subgroup for an element in a flat group with larger flat rank, it does not contain elements of the larger flat group. 
\begin{example}
Let	$A := \langle a,b\rangle$ be the free group on generators $a$ and $b$ and define a homomorphism $\phi : A\to \Aut(\mathbf{Q}_p^2)$ by 
$$
\phi(a)(\xi_1,\xi_2) = (p\xi_1,\xi_2)\text{ and }\phi(b)(\xi_1,\xi_2) = (\xi_1,p\xi_2),\qquad (\xi_1,\xi_2)\in\mathbf{Q}_p^2.
$$
Then $A$ is a flat subgroup of $G := \mathbf{Q}_p^2\rtimes A$ having $\mathbf{Z}_p^2\rtimes \{0\}$ as tidy subgroup. The uniscalar subgroup of $A$ is the commutator subgroup of $\langle a,b\rangle$ and the flat rank of $A$ is equal to $2$.

Let $H = \langle a\rangle\leq A$. Then the flat rank of $H$ is equal to $1$. However, while we have $\lev{H} = \{0\}\times \mathbf{Q}_p = \con{b}$ and $s(b^{-1}|_{\lev{H}}) = p$, the set of $a$-conjugates of $b^{-1} $ is infinite and discrete, and hence does not have compact closure. Hence $b^{-1}$ does not belong to $\lev{H}$. 
\end{example}

\section{Precompact conjugation orbits and a symmetric relation}  
\label{sec:symmetric_relation}

The relation \emph{commutes with}, as in `$x_1$ commutes with $x_2$', is reflexive and symmetric on the elements of a group. The weaker relation \emph{in the Levi subgroup of}, as in `$x_1\in \lev{x_2}$', is reflexive but is not symmetric. Theorem~\ref{thm:extend_flat} and ideas in its proof suggest a still weaker reflexive relation that is symmetric. That relation, called \emph{vague commutation} and denoted $x_1Vx_2$, is explored in this section. Unfortunately, symmetry of the relation is achieved at the expense of $\{x_1\in G \mid x_1V x_2\}$ being a group.
\begin{definition}
	Let $G$ be a \tdlc~group. For $x_1,x_2\in G$ say that \emph{$x_1$ vaguely commutes with $x_2$}, and write $x_1Vx_2$, if there are $n>0$ and a compact subgroup, $U$, of $G$ normalised by $x_1^n$ and $u\in U$ such that $ux_1^{n}\in\lev{x_2}$.  
\end{definition}

It is clear that the relation $VC$ is reflexive and the next proposition shows that it is symmetric as well.
\begin{proposition}
	\label{prop:relation_is_symmetric}
	Suppose that $x_1,x_2\in G$ and that $x_1Vx_2$. Then there are, for each $i=1,2$, integers $n_i>0$, compact subgroups, $U_i$, normalised by $x_i^{n_i}$  and elements $u_i\in U_i$, such that $z_i := u_ix_i^{n_i}$ commute modulo a compact subgroup, $L$, satisfying that $z_iLz_i^{-1} = L$. 
	
	Moreover, $x_2Vx_1$. 
\end{proposition}
\begin{proof}
Since $x_1$ vaguely commutes with $x_2$, there are $n_1>0$ and a compact  subgroup, $U_1$ of $G$, normalised $x_1$ and $u_1\in U_1$ such that $u_1x_1^{n_1}\in\lev{x_2}$. Then $H := \langle x_2\rangle$ and $y := u_1x_1^{n_1}$ satisfy the hypotheses of Theorem~\ref{thm:extend_flat}. Hence there are $n_2>0$ and a compact subgroup, $U_2$, normalised by $x_2$ and $u_2\in U_2$, and a compact subgroup, $L$, such that, setting $z_i = u_ix_i^{n_i}$: $z_1Lz_1^{-1} = L = z_2Lz_2^{-1}$ and $[z_1,z_2] \in L$, as claimed. 

It follows that $u_2x_2^{n_2}\in\lev{u_1x_1^{n_1}}$ and thus that
$$
\left\{(u_1x_1^{n_1})^k (u_2x_2^{n_2})(u_1x_1^{n_1})^{-k}\mid k\in\mathbf{Z}\right\}
$$
has compact closure. Since $u_1\in U_1$, which is normalised by $x_1^{n_1}$, we have that $(u_1x_1^{n_1})^k = \tilde{u}_kx_1^{kn_1}$ for each $k$ with $\tilde{u}_k\in U_1$. Since $U_1$ is compact, it follows that
$$
\left\{x_1^{kn_1} (u_2x_2^{n_2}) x_1^{-kn_1}\mid k\in\mathbf{Z}\right\}
$$
has compact closure, and hence so does $\left\{x_1^{k} (u_2x_2^{n_2}) x_1^{-k}\mid k\in\mathbf{Z}\right\}$. Therefore $u_2x_2^{n_2}$ belongs to $\lev{x_1}$ and $x_2$ vaguely commutes with $x_1$. 
	\end{proof}

The \emph{vague centraliser} of $x_2$ is 
$$
\left\{ x_1\in G \mid x_1\text{ vaguely commutes with }x_2\right\}.
$$  
Then the vague centraliser of $x_2$ is the analogue of the centraliser of $x_2$ in the case of the relation of commutation, and of $\lev{x_2}$ for the `belongs to $\lev{x_2}$' relation. 

\begin{example}
	\label{examp:vague_centraliser}
Let $\Aut(\tree_{q+1})$ be the automorphism group of a regular tree $\tree_{q+1}$ with valency $q+1$ and $q\geq2$. Let $x_1$ and $x_2$ be translations, that is, automorphisms having infinite orbits, with axes $\ell_1$ and $\ell_2$, see {\color{red}\cite[]{Tits_arbre}}. Then $x_1^kx_2x_1^{-k}$ is a translation with axis $x_1^k.\ell_2$, and the same is true if $x_2$ is replaced by $ux_2^n$ for any $n>0$ and $u$ in a compact subgroup normalised by $x_2^n$. Hence $\{x_1^k(ux_2^n)x_1^{-k}\}_{k\in\mathbf{Z}}$ does not have compact closure unless $\ell_1 = \ell_2$, and it follows that $x_1$ and $x_2$ do not vaguely commute unless this holds. On the other hand, if $x_1$ is an elliptic element of $\Aut(\tree_{q+1})$, that is, if $x_1$ has finite orbits in $\tree_{q+1}$, then $\langle x_1\rangle^-$ is compact. Hence so is $\{x_1^kx_2x_1^{-k}\}^-_{k\in\mathbf{Z}}$,	
and $x_1$ and $x_2$ vaguely commute. 

These considerations thus show that the vague centraliser of $x$ is: 
\begin{itemize}
\item $\Aut(\tree_{q+1})$, if $x$ is elliptic; and 
\item $\{\text{translations with axis }\ell\} \cup \{\text{elliptic elements}\}$, if $x$ is a translation with axis $\ell$.
\end{itemize}
\end{example}

Example \ref{examp:vague_centraliser} shows that, while the centraliser of $x_2$ and $\lev{x_2}$ are subgroups of~$G$, the vague centraliser is not in general. In order to expand on the example, we make a few more definitions.
\begin{definition}
	\label{defn:Centre-like_sets}
Let $G$ be a totally disconnected, locally compact group. 
\begin{enumerate}
	\item  ${Per}(G) = \{x\in G\mid \langle x\rangle^-\text{ is compact}\}$.
	\item $x\in G$ is an \emph{${FC}^-$-element} if the conjugacy class of $x$ has compact closure, and the set of ${FC}^-$-elements in $G$ is denoted ${FC}^-(G)$.
	\item The \emph{approximate centraliser of $G$} is the set of all $x\in G$ such that $\lev{x} = G$ and is denoted $AC(G)$.
	\item The \emph{approximate centre of $G$} is equal to $\bigcap\left\{\lev{y}\mid y\in G\right\}$ and is denoted $AZ(G)$. 
\item The \emph{vague centre} of $G$ is the set of all $x\in G$ whose vague centraliser is equal to $G$, and is denoted $VZ(G)$.
\end{enumerate}
\end{definition}

In Example \ref{examp:vague_centraliser}, ${Per}(G)$ is equal to the set of elliptic automorphisms of $\tree_{q+1}$, and this set is equal to the approximate centraliser $AC(G)$ and to the vague centre $VZ(G)$. The set of ${FC}^-$-elements ${FC}^-(G)$ and the approximate centre $AZ(G)$ are equal to $\triv$. 

The next result states the order relationships between the sets in Definition~\ref{defn:Centre-like_sets} and records their properties. 
\begin{figure}[h]
	\label{fig:Centres}
 \begin{tikzpicture}
\path (-1.5,-1.5)   node(A) [rectangle]  {${Per}(G)$}
         (1.5,1.5)   node(B)  [rectangle]  {$AC(G)$}
         (3,-1.5)  node(C)  [rectangle]  {$Z(G)$}
         (4.5,0)  node(D)  [rectangle]  {${FC}^-(G)$}
         (3,3)  node(F)   [rectangle]  {$VZ(G)$}
         (4.5,1.5)  node(E)    [rectangle]  {$AZ(G)$};
         
         \draw[semithick] (A) -- (B) -- (F);
         \draw[semithick] (B) -- (C) -- (D) -- (E) -- (F);
      \end{tikzpicture}
\caption{Subsets of $G$ related to the centre}
\end{figure}
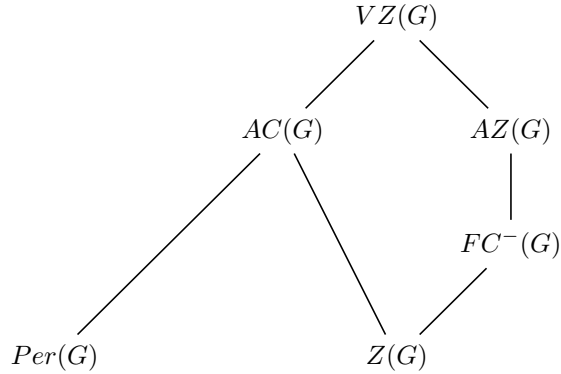

\begin{proposition}
	\label{prop:vague_centraliser}
Let $G$ be a totally disconnected, locally compact group. The subsets of a totally disconnected, locally group, $G$, described in Definition~\ref{defn:Centre-like_sets} are ordered by inclusion as shown the Hasse diagram in Figure~\ref{fig:Centres}. Moreover: 
\begin{enumerate}
	\item all subsets are invariant under conjugation;
	\item $Z(G)$, ${FC}^-(G)$ (provided that $G$ is compactly generated) and $AZ(G)$ are closed normal subgroups of $G$; and
	\item $\text{Per}(G)$, $AC(G)$ and $VZ(G)$ are closed subsets of $G$.
\end{enumerate}
The scales of all elements in $VZ(G)$ are equal to $1$. 
\end{proposition}
\begin{proof} The inclusions may be seen as follows: Proposition~\ref{prop:relation_is_symmetric} shows that $AC(G)$ and $AZ(G)$ are contained in $VZ(G)$; $AC(G)$ contains $Per(G)$ because compactness of $\langle x\rangle^-$ implies that $\{x^kyx^{-k}\}^-_{k\in\mathbf{Z}}$ is equal to the set of conjugates of $y$ by $\langle x\rangle^-$, and contains $Z(G)$ because $\{x^kyx^{-k}\}^-_{k\in\mathbf{Z}} = \{y\}$ if $x\in Z(G)$; $AZ(G)$ contains $FC^-(G)$ because $\{y^kxy^{-k}\}^-_{k\in\mathbf{Z}}\subseteq \{wxw^{-1}\}^-_{w\in G}$ for every $y\in G$; and $FC^-(G)$ contains $Z(G)$ because singleton conjugacy classes are compact.
	
Since each of the statements defining the sets in Definition~\ref{defn:Centre-like_sets} is invariant under conjugation, the sets so defined are too. 
	
For the further properties of the sets, it is well-known that the centre is a closed subgroup in every topological group, and that $FC^-(G)$ is closed if $G$ is compactly generated is \cite[Theorem~2]{Moller_FC}. That $AZ(G)$ is a closed subgroup holds because it is the intersection of the closed subgroups $\lev{y}$, $y\in G$. That these subgroups are normal holds because, as just seen, they are closed under conjugation. That ${Per}(G)$ is closed is \cite[Theorem~2]{Wi:tdHM}. The proofs that $AC(G)$ and $VZ(G)$ are closed too mirror the argument used for $Per(G)$, which uses the fact that the group elements are uniscalar. 

Assume, as will shortly be proved, that all elements of $x\in VZ(G)$ are uniscalar and consider $y\in VZ(G)^-$. To see that $y\in VZ(G)$, it must be shown that for every $y_2\in G$ there are $n>0$, a compact subgroup, $U$, normalised by $y_2^n$ and $u\in U$ such that $uy_2^n\in\lev{y}$. We have that $s(y) = 1 = s(y^{-1})$, because the scale is continuous, and hence there is a compact, open $W\leq G$ normalised by $y$. Since $yW$ is a neighbourhood of $y$, we may choose $x\in VZ(G)\cap yW$. Then there are $n>0$, a compact subgroup, $U$, normalised by $y_2$ and $u\in U$ such that $uy_2^n\in\lev{x}$ because $x\in VZ(G)$. Since $yW = xW$ and $x$ normalises $W$, there are elements $w_k\in W$, for all $k\in\mathbf{Z}$, such that $y^k = w_kx^k$. Hence, for every $k\in\mathbf{Z}$,
$$
y^k(uy_2^n)y^{-k} = w_k x^k(uy_2^n)x^{-k}w_k^{-1}\in W\{x^j(uy_2^n)x^{-j}\}_{j\in\mathbf{Z}}W,
$$
which is compact. Hence $uy_2^n\in \lev{y}$ as desired and $y\in VZ(G)$. The proof that $AC(G)$ is closed is essentially the same. Every element of $AC(G)$ is uniscalar because it is contained in $VZ(G)$. Given $y\in AC(G$, there is $W$ normalised by $y$ and $x\in AC(G)\cap yW$. Then $\lev{x}= \lev{y}$ and the fact that every $y_2\in G$ belongs to $\lev{x}$ implies that $y_2$ also belongs to $\lev{y}$. 

The proof will now be completed by showing that every element of $VZ(G)$ is uniscalar. Let $x_1$ in $G$ have $s(x_1)>1$, choose $v_+\in\con{x}\setminus\nub{x}$, and put $x_2 = v_+x_1$.  Then for every  $n>0$, every compact group $U$ normalised by $x_1^n$ and $u\in U$, and for all $k\geq0$, we have
\begin{align*}
&\phantom{=}\ x_2^{k+1}(ux_1^n)x_2^{-k-1} = (v_+x)^{k+1} (ux_1^n)(v_+x)^{-k-1} \\
&= v_+\dots (x_1^kv_+x_1^{-k})(x_1^{k+1}ux_1^{-k-1})(x_1^{n+k}v_+x_1^{-n-k})^{-1}\dots (x_1^nv_+x_1^{-n})^{-1}x^n.
\end{align*}
Hence the sequence $\{x_2^{k+1} (ux_1^n)x_1^{-k-1}\}_{k\geq0}$ is discrete and non-compact, by Lemma~\ref{lem:unbounded_sequence}. In particular, $ux_1^n$ does not belong to $\lev{x_2}$ for any $n$ and $u$, and $x_1$ and $x_2$ do not vaguely commute. Since for any $x_1$ with $s(x_1)>1$ there is $x_2$ that does not vaguely commute with $x_1$, every element in $VZ(G)$ is uniscalar. 
\end{proof}  

\begin{remark}
(1) If $x_1\in \text{Per}(G)$, then $x_2\in \lev{x_1}$, and hence $x_2$ vaguely commutes with $x_1$, for every $x_2$ in $G$. Proposition~\ref{prop:relation_is_symmetric} then implies that every $x_1\in \text{Per}(G)$ vaguely commutes with every $x_2\in G$. It is instructive to see this directly: consider $x_1\in \text{Per}(G)$ and $x_2\in G$, and take $n=1$, $U =  \langle x_1\rangle^-$ and $u=x_1^{-1}$; then $ux_1 = \ident_G$, which certainly belongs to $\lev{x_2}$. \\
(2) The subgroup $FC^-(G)$ need not be closed if $G$ is a \tdlc~group that is not compactly generated. The example in \cite[Proposition~3]{Tits_FC} is a group $G$ in which $FC^-(G)$ is proper and dense. If $G$ is a connected locally compact group, then $FC^-(G)$ is closed, by \cite[Corollary~1]{Tits_FC}, and more information about the structure of this subgroup is given there. The structure of general $[FC]^-$-groups is described in \cite[Theorem 3D]{Robertson_FCbar}, see \cite[Theorem 2.2]{Liukkonen_FCbar} for a proof and also \cite[12.6.10]{PalmerII}.\\
(3) The proofs that $VZ(G)$ and $AC(G)$ are closed are modelled on the argument used in \cite[Theorem~2]{Wi:tdHM} to show that $Per(G)$ is closed, and that theorem answers a question posed by Karl Heinrich Hofmann.\\
(4) The class of all compactly generated, topologically simple, non-discrete locally compact groups is denoted by $\mathcal{S}$ and it is an open question whether there any \tdlc~group in $\mathcal{S}$ can be uniscalar. It follows from Proposition~\ref{prop:vague_centraliser} that a sharper version of this question is whether any groups in $\mathcal{S}$ have $G$ equal to $VZ(G)$, $AC(G)$ or $Per(G)$. Any group $G$ equal to $Per(G)$ has the stronger property than uniscalarity of being \emph{anisotropic}, which means all contraction subgroups for inner automorphisms are trivial. Note that discrete finitely generated, simple groups with $G=Per(G)$  exist because, for discrete groups, $Per(G)$ is the set of torsion elements and there are simple counterexamples to Burnside's conjecture. 
\end{remark}

\bibliographystyle{plain}

\end{document}